\documentclass{amsart}
\usepackage{epsfig}
\usepackage{latexsym}
\usepackage{amscd}
\usepackage{amssymb}
\usepackage{amsthm}
\usepackage{amsmath}
\usepackage{amsxtra}
\usepackage{amsfonts}
\usepackage{mathrsfs}

\usepackage[all]{xy}

\newtheorem{prop}{Proposition}
\newtheorem{thm}{Theorem}
\newtheorem{cor}{Corollary}

\newtheorem{lemma}{Lemma}

\theoremstyle{remark}

\newtheorem{exm}{Example}
\DeclareMathOperator{\Frob}{Frob}

\DeclareMathOperator{\Ext}{Ext}

\def\be{\kern -.1em}
\def\lbe{\kern -.025em}

\DeclareFontEncoding{OT2}{}{} 
\newcommand{\textcyr}[1]{%
 {\fontencoding{OT2}\fontfamily{wncyr}\fontseries{m}\fontshape{n}
 \selectfont #1}}
\newcommand{\Sha}{{\!\be\lbe\mbox{\textcyr{Sh}}}}

\newcommand{\A}{A}
\newcommand{\B}{B}
\newcommand{\E}{E}

\newcommand{\G}{\mathbb G}
\newcommand{\T}{T}

\newcommand{\Proj}{\mathbb P}
\newcommand{\Z}{Z}

\newcommand{\X}{X}
\newcommand{\Y}{Y}

\newcommand{\Gm}{\G_{\text m}}

\newcommand{\la}{\langle}
\newcommand{\ra}{\rangle}
\newcommand{\rarr}{\rightarrow}

{.tfm}
{.tfm}

\title[Product of local points of subvarieties of semi-abelian varieties]{Product of local points of subvarieties of almost isotrivial semi-abelian varieties over a global function field}

\author{Chia-Liang Sun}

\begin{document}
\maketitle

\begin{itemize}
\item Physical Address: Institute of Mathematics, Academia Sinica,
6F, Astronomy-Mathematics Building,
No. 1, Sec. 4, Roosevelt Road,
Taipei 10617, TAIWAN
\item E-mail Address: {\tt csun@math.sinica.edu.tw}
\end{itemize}

\begin{abstract}
For a semi-abelian variety over a global function field which is isogenous to an isotrivial one, we show that on the product of local points of
a subvariety satisfying a minor condition, 
the topological closure of a finitely
generated subgroup of 
global points cuts out
exactly the global points of the subvariety lying in this subgroup. As a corollary, on every non-isotrivial super-singular curve of genus two over a global function field, we conclude that the Brauer-Manin condition cuts out exactly the set of its rational points.
\end{abstract}

\section{Introduction}\label{AdelicIsotrivial_intro}
Let $K$ be a global function field of characteristic $p$, that is, a finitely generated field extension over the prime field $\mathbb{F}_p$ with the
transcendence degree $1$.
We fix an algebraic closure $\overline{K}$ and let $K^{s}$ denote the separable closure of $K$.
Let $\Omega_K$ denote the set of all places of $K$ and let $\Omega$ denote a co-finite subset of $\Omega_K$.
For for each $v\in \Omega_K$, let $K_v$ denote the completion of $K$ at $v$.
For an algebraic variety $X$ defined over $K$, we endow $X(K_v)$ with the natural $v$-adic topology,
and then endow $\prod_{v\in\Omega}X(K_v)$
with the product topology.
In this paper, we assume that $X$ is a closed subvariety of a given semi-abelian variety $\A$, both  are defined over $K$.
We identify every subset $H\subset \A(K)$ as a topological subspace of $\prod_{v\in\Omega}\A(K_v)$ by the diagonal embedding, and denote by $\overline{H}$ its topological closure; moreover, for each $v\in\Omega_K$,
the inclusion $H\rarr A(K_v)$ is continuous and therefore induces a subtopology of $H$, which  will be referred to as \em $v$-adic subtopology\em.

Suppose $H$ is a subgroup of $A(K)$. The main aim of this paper is to investigate the circumstances where the equality
\begin{equation}\label{maineq}
\prod_{v\in\Omega}\X(K_v)\cap\overline{H}=\X(K)\cap H
\end{equation}
holds. To do so, we shall assume that $H$ is {\em{finitely generated}}, as Example \ref{exm0} in Subsection \ref{pfbase} shows that
the equality would not hold even in the simplest case where $\A=\mathbb{G}_m$, otherwise.
On the other hand, the case where $A$ is an abelien variety and $H=A(K)$,
has been studied by Poonen and Voloch \cite{PoonenVoloch}. Indeed, they propose that
in general the equality
$$
\prod_{v\in\Omega}\X(K_v)\cap\overline{\A(K)}=\overline{\X(K)}
$$
should hold, and they prove \eqref{maineq} under the hypothesis that
$A_{\overline{K}}$ has no isotrivial quotient, $A(K^s)[p^{\infty}]$ is finite and $X$ is coset-free
({\em{op. cit.}}, Conjecture C and Theorem B).
In this paper, we consider the case where $A$ is  isogenous to an isotrivial semi-abelian variety.
Recall that in the category of varieties (resp. of semi-abelian varieties), an object defined over $K$ is \em isotrivial \em if it is isomorphic over $\overline{K}$ to one defined over $\overline{\mathbb{F}_p}$.
Our main result is the following:

\begin{thm}\label{AdelicIsotrivial_main} Let $X$ be a closed subvariety of a semi-abelian variety $\A$, both  are defined over $K$. Assume that there is an isogeny $f$ defined over $\overline{K}$
from $\A$ to a semi-abelian variety $\A_0$ defined over $\overline{\mathbb{F}_p}$
so that each translate $P+f(\X)$, $P\in \A_0(\overline{K})$, of $f(X)$ contains no positive-dimensional closed subvariety which is $\overline{\mathbb{F}_p}$-rational  in $(\A_0)_{\overline{\mathbb{F}_p}}$.
Then for every finitely generated subgroup $H$ of $\A(K)$, the set $\X(K)\cap H$ is finite and the equality \eqref{maineq} holds.
\end{thm}
Example \ref{exm1} at the end of Section \ref{ind} explains why we cannot expect the conclusion of our main result to
hold without any assumption on $\X$.
The proof of the theorem consists of two parts which are carried out respectively in Section \ref{base} and Section \ref{ind},
with the key ingredients, Proposition \ref{CSPmain}, Lemma \ref{int} and Lemma \ref{gen}, proved in Section \ref{pfkey}.
The first part treats the case
where $\X$ is zero-dimensional by adapting
the proof of Proposition 3.7 in \cite{PoonenVoloch} to our situation, while the second reduces the general case to the
zero-dimensional case by the induction on the dimension of $\X$ using a different Mordell-Lang type argument.
The approach for the second part is originated from the proof of
Theorem A. Part 1 in \cite{AV92}.

Suppose $\A$ is the Jacobian of $\X$ which is embedded into $\A$ under an Albanese map induced from a divisor on $\A$ defined over $K$ of degree $1$.
It is proved in \cite{PoonenVoloch}, Section 4, that if the Tate-Shafarevich group $\Sha(\A)$ of $\A$ is finite,
then the set
$\prod_{v\in\Omega_K}\X(K_v)\cap\overline{\A(K)}$ 
is in bijection with the Brauer set of $\X$ over $K$.
Therefore, if Theorem \ref{AdelicIsotrivial_main} holds, then the Brauer-Manin condition cuts exactly the set of its $K$-rational points on $\X$.
An non-isotrivial projective curve can have its Jacobian isogenous to an isotrivial abelian variety; for examples, super-singular curves of genus $2$ have this property \cite{CurVar}. The following corollary is proved in Section \ref{ind}.

\begin{cor}\label{ss2curve}
On any non-isotrivial projective $K$-curve with its Jacobian isogenous to an isotrivial abelian variety, the Brauer-Manin condition cuts exactly the set of its $K$-rational points.
\end{cor}

\section{Preliminaries}\label{pre}
\subsection{}\label{quasi_proj}

A semi-abelian variety $\A$ is said to be defined over a subfield $F$ of $\overline{K}$ if the  underlying variety has the following form
\begin{equation}\label{qvar}
\left\{P\in\Proj^N(\overline{K}):
\begin{array}{l}
f(P)=0\text{ for all }f\in I\\ g(P)\neq 0\text{ for some }g\in J
\end{array}
\right\},
\end{equation}
where $I$ and $J$ are homogeneous ideals in the polynomial ring $F[X_0,\ldots,X_N]$, such that the group laws are expressible in $F$.
For any field $L$ such that $F\subset L \subset \overline{K}$, we say a closed subvariety of $\A$ is \em $L$-rational in $\A_F$ \em if, whenever $\A$ is expressed as (\ref{qvar}), it has the form
\begin{equation*}
\left\{P\in\Proj^N(\overline{K}):
\begin{array}{l}
f(P)=0\text{ for all }f\in \tilde{I}\\ g(P)\neq 0\text{ for some }g\in J
\end{array}
\right\}
\end{equation*}
for some homogeneous ideal $\tilde{I}$ in $L[X_0,\ldots,X_N]$.

For each $v\in \Omega_K$, we denote by $O_v$
the valuation ring in $K_v$ and by $m_v$ the
maximal ideal in $O_v$.
For a finite subset $S$ of $\Omega_K$, the subring of $S$-integers in $K$ is denoted by
$$
O_S=\left\{x\in K: x\in O_v \text{ for each } v\notin S \right\}.
$$
Given a such $\A$ defined over $K$ as above, let $J^{(v)}$ denote the subset of $J\cap O_v[X_0,\ldots,X_N]$ consisting of polynomials with some coefficients not in $m_v$,
and define
\begin{equation}\label{VOv}
\A(O_v)=\left\{P\in\Proj^N(K_v):
\begin{array}{l}
P=[x_0:\ldots:x_N]\\
x_i\in O_v \text{ for all } i\\
x_{i_0}\notin m_v \text{ for some } i_0\\
f(x_0,\ldots,x_N)=0 \text{ for all }f\in I\\
g(x_0,\ldots,x_N)\notin m_v \text{ for some }g\in J^{(v)}
\end{array}
\right\}.
\end{equation}

Also, for any finite subset $S\subset\Omega$, set
\begin{equation}\label{VOS}
\A(O_S)=\bigcap_{v\in \Omega_K\setminus S}\left(A(K)\cap A(O_v)\right).
\end{equation}
Then, the following useful result can be obtained by Hilbert's Nullstellensatz:
\begin{lemma}\label{preOv}
Let $\phi: \A\rarr\A$ be a regular map defined over $K$. 
Then $\phi$ preserves $\A(O_v)$ for all but finitely many $v\in\Omega_K$, and hence $\A(O_S)$ for $S$ sufficiently large.\qed
\end{lemma}


\subsection{}\label{topology} The group operations on $\A$, given by regular maps defined over $K$, are continuous with respect to the topology on $\A(K_v)$, for each $v\in\Omega_K$. This implies that $\A(K_v)$ is a Hausdorff topological abelian group.
It is totally disconnected since $v$ is non-archimedean. Also, (\ref{VOv}) shows that $\A(K_v)$ contains
$\A(O_v)$ as a compact open subset. Consequently, it is locally compact, and hence complete.
By Theorem (7.7), \cite{HarAna}, it follows that the topology of $\A(K_v)$ is generated by open subgroups, and therefore so are $\prod_{v\in\Omega}\A(K_v)$ and its subgroups.

\begin{lemma}\label{adehaus}
Every finitely generated subgroup $H$ of $\A(K)$ admits a Hausdorff subtopology generated by subgroups with finite index.
\end{lemma}
\begin{proof}
Because $H$ is finitely generated and $\Omega_K\setminus\Omega$ is finite, there exists a place $v_0\in\Omega$ such that $\A(O_{v_0})$ is a group containing $H$. Since $\A(O_{v_0})$ is a compact subgroup of $\A(K_{v_0})$ whose topology is Hausdorff and generated by subgroups, the topology of $\A(O_{v_0})$ is  Hausdorff and  generated by subgroups with finite index, and so is the $v_0$-adic subtopology of $H$.
\end{proof}

Let $\mathfrak{P}(A,X,f,K)$ stand for the statement of Theorem \ref{AdelicIsotrivial_main}.

\begin{lemma}\label{1stred} Let $A,X,f,K$ be as in {\em{Theorem \ref{AdelicIsotrivial_main}}} and let $L/K$ be a finite
extension. Then $\mathfrak{P}(A,X,f,L)\Rightarrow\mathfrak{P}(A,X,f,K)$.
\end{lemma}
\begin{proof} Note that the definition of $\overline{H}$ depends on the choice of $K$.
Let $\Omega_L^0$ denote the set of places of $L$ lying above $\Omega$.
Consider the natural embeddings
$$\xymatrix{ H \ar@{^{(}->}[r]  &  \prod_{v\in\Omega}\A(K_v)\ar@{^{(}->}[r]^i & \prod_{u\in\Omega_L^0}\A(L_u)}.$$
Since $i$ actually identifies $\prod_{v\in\Omega}\A(K_v)$ as a closed subgroup of $\prod_{u\in\Omega_L^0}\A(L_u)$,
$\overline{H}$ will be remained the same, if $K$ is replaced by $L$.
If $\mathfrak{P}(A,X,f,L)$ holds, then
$\prod_{w\in\Omega_L^0}\X(L_w)\cap\overline{H}=\X(L)\cap H$ is a finite set. This implies
\begin{equation*}
\begin{array}{cl}
 &\prod_{v\in\Omega}\X(K_v)\cap\overline{H}\\
=&\left(\prod_{w\in\Omega_L^0}\X(L_w)\cap\prod_{v\in\Omega} \A(K_v)\right)\cap\overline{H}\\
=&\prod_{w\in\Omega_L^0}\X(L_w)\cap\overline{H}\\
=&\X(L)\cap H\\
=&\left(X(L)\cap\A(K)\right)\cap H\\
=&\X(K)\cap H,
\end{array}
\end{equation*}
which is also a finite set.
\end{proof}
\subsection{}\label{semi_abel} In view of Lemma \ref{1stred}, by replacing $K$ by certain finite extension $L$ if necessary,
we can assume that $A$ is an extension of an abelian variety $B$ by a split torus $\Gm^n$, for some $n\geq 0$. This is actually the
definition taken by Serre \cite{AGCF}. Lemma \ref{preOv} implies that $\A(O_v)$ is a group for all but finitely many $v\in\Omega_K$.
Therefore, if $S$ is sufficiently large, then $\A(O_S)$ is a subgroup of $\A(K)$. In this case, it is finitely generated.
To see this, we may extend $S$ and assume that the subgroup $\Gm^n(O_S)\subset A(O_S)$ coincides with $(O_S^*)^n$.
Since $\A(O_S)$ is mapped into $\B(K)$ with $\Gm^n(O_S)$ as the kernel, the assertion follows from
the Mordell-Weil theorem and Dirichlet's unit theorem.
Also, if $S$ is large enough, then the finitely generated subgroup $H$ is contained in $\A(O_S)$.




\begin{lemma}\label{CSPuseful}
Every isotrivial semi-abelian variety $A_0$ is isogenous {\em{(over $\overline{K}$)}}
to the product $\Gm^{n_0}\times \B_0$, for some nonnegative integer $n_0$, where $\B_0$ is an abelian variety
defined over $\overline{\mathbb{F}_p}$.
\end{lemma}
\begin{proof} We may assume that $A_0$ is actually defined over $\overline{\mathbb{F}_p}$. Then there is a strictly exact sequence (see \cite{AGCF})
$1\rarr \Gm^{n_0}\rarr \A_0 \rarr \B_0\rarr 0$, in which $\B_0$ is an abelian variety defined over $\overline{\mathbb{F}_p}$.
It remains to show that
$\A_0$ is isogenous to $\Gm^{n_0}\times\B_0$.

Recall that, for any pair $(\G,\T)$ of commutative algebraic groups, the set of isomorphism classes of
commutative algebraic groups $\E$ along with the strictly exact
sequence $0\rarr \T\rarr \E \rarr \G\rarr 0$ forms an abelian group
$\Ext(\G,\T)$ under the Baer sum. In fact, $\Ext$ is a bifunctor
from the category of pairs of commutative algebraic groups to the
category of abelian groups \cite{AGCF}. It is routine to check that, for any
commutative algebraic group $\E$ representing its class
$[\E]\in\Ext(\G,\T)$, and any positive integer $m$, there is a
natural exact sequence
\begin{equation}\label{isoexact}
0\rarr \T[m]\rarr \E \rarr \E^{(m)}\rarr 0
\end{equation}
defined over an algebraic closure of a field of definition of $\E$, where $\T[m]$ is the algebraic subgroup of $m$-torsion points in $\T$, and $\E^{(m)}$ is a commutative algebraic group representing the class $m[\E]\in\Ext(\G,\T)$. In
our case, the isomorphism class $[\A_0]$ of $\A_0$ lies in $\Ext(\B_0,\Gm^{n_0})=\Ext(\B_0,\Gm)^{n_0}$, where the equality holds because $\Ext$ is a bifunctor. One knows (e.g., the comments following Theorem 6 of Chapter VII in \cite{AGCF}) that $\Ext(\B_0,\Gm)$
is isomorphic to the dual abelian variety $\B_0'$ of $\B_0$. In particular, $[\A_0]$ lies in $(\B_0'(\overline{\mathbb{F}_p}))^{n_0}$
which is a torsion group. Therefore $m[\A_0]=0$, for some $m$, and hence $\Gm^{n_0}\times\B_0=A_0^{(m)}$ and the isogeny is provided by
\eqref{isoexact}.
\end{proof}

\section{The Zero-Dimensional Case}\label{base}
In this section, we prove Theorem \ref{AdelicIsotrivial_main}, by assuming $\dim\X=0$.

\subsection{A Uniform Filtration Over All $v$-adic subtopologies}\label{unifil}
Suppose $A_0$ is a semi-abelian variety define over a finite field $\mathbb{F}_q$ containing $\mathbb{F}_p$.
Then the Frobenius morphism $\Frob:\A_0\rarr\A_0$ is well defined, and if $\A_0$ is embedded as a subvariety of $\Proj^N$
then $\Frob$ is simply the restriction of the Frobenius map  on $\Proj^N$,
sending $[x_0:\ldots:x_N]$ to $[x_0^q:\ldots:x_N^q]$.
Thus, $\Frob$ preserves the group structure on $\A_0$ and it induces an injective map $\A_0(O_S)\rarr\A_0(O_S)$, which we also denote by $\Frob$
which is a group homomorphism, if $\A_0(O_S)$ is a group.

\begin{prop}\label{frob}
Suppose that $f:\A\rarr \A_0$ is an isogenous defined over $\overline{K}$ and $A_0$ is an isotrival semi-abelian variety. Then for any finitely generated subgroup $H$ of $\A(K)$, there exists a collection $\{U_n:n\geq 1\}$ of subgroups of $H$ with the following two properties:
\begin{enumerate}
    \item\label{uniopen} For each $v\in\Omega_K$ and each $n\geq 1$, $U_n$ is open in every $v$-adic subtopology of $H$.
    \item\label{inttor} $\bigcap_{n\geq 1} U_n$ is contained in the torsion subgroup of $H$.
\end{enumerate}
\end{prop}
\begin{proof}
Without loss of generality, we may replace $K$ by one of its finite extensions and assume that $f$ is defined over $K$ and $A_0$ is defined over
$\mathbb{F}_q\subset K$. Since $\xymatrix{H\ar[r]^f & f(H)}$ is continuous with finite kernel, it is sufficient to show that $f(H)$ precesses
a family of open sets satisfying the properties correspondingly. Thus, we only need to consider the case where $\A=\A_0$ and $f$ is the identity map.

We claim that for each $v\in\Omega_K$ and each $n\geq 1$, the subgroup $U_n:=\Frob^n(\A(K))\cap H\subset H$ is open in the $v$-adic subtopology.
Then, (\ref{inttor}) holds, as
$$\bigcap_{n\geq 1} U_n\subset \bigcap_{n\geq 1}\Frob^n(\A(K))\subset \A(\bigcap_{n\geq 1}K^{p^n})$$
and $\bigcap_{n\geq 1}K^{p^n}$ is the maximal finite subfield of $K$.

To prove the claim, we first note that since there is no nontrivial purely inseparable finite extension of $K$ inside $K_v$,
$\Frob^n(\A(K_v))\cap \A(K)\subset \Frob^n(\A(K))$, and hence $\Frob^n(\A(K_v))\cap H\leq U_n$.
Then it remains to show that $\Frob^n(\A(K_v))\cap H$ is open in the $v$-adic subtopology of $H$. It is clear that that $\Frob^n(\A(K_v))$ is closed in $\A(K_v)$, and consequently the quotient space $\A(K_v)/\Frob^n(\A(K_v))$ is Hausdorff. Consider the map $H\rarr \A(K_v)/\Frob^n(\A(K_v))$ induced from the inclusion $H\subset\A(K_v)$. It is continuous with respect to $v$-adic subtopology of $H$.
Also, as it factors through $\A(O_S)/\Frob^n(\A(O_S))$, for some finite $S\subset \Omega_K$ such that $H\leq \A(O_S)$,
the image of the map is finite, whence discrete. This completes our proof.
\end{proof}

\subsection{Congruence Subgroup Property}\label{CSP}
For an additive topological abelian group $G$, we say that $G$ has the \em {\bf congruence subgroup property} \em if the subgroup $nG=:\{nP: P\in G\}$ is open for every positive integer $n$; if $G$ is finitely generated, then the following  conditions are equivalent:

\begin{enumerate}
    \item $G$ has the congruence subgroup property.
    \item Every subgroup of $G$ of finite index is open.
    \item Every subgroup of $G$ is closed.
\end{enumerate}

\begin{lemma}\label{CSPsub} Let $G$ be a finitely generated abelian topological groups. Let $\Sigma$ be a set consisting of natural numbers,
which is closed under multiplication, satisfying the condition that every subgroup of $G$ of index in $\Sigma$ is open.
Then $\Sigma$ also satisfying the corresponding condition for each subgroup $H$ of $G$, namely, every subgroup of $H$ of index in $\Sigma$
is open in $H$. In particular, if $G$ has the congruence subgroup property, then so has $H$.
\end{lemma}
\begin{proof} 
We may assume that each positive divisor of any element in $\Sigma$ also lies  in $\Sigma$. Let $m\in\Sigma$ be the product of those natural numbers in $\Sigma$, each of which is the order of some elements in the finitely generated abelian group $G/H$. Then, for every  $n\in\Sigma$, we have $H\cap mnG\leq nH$. Since $mn\in\Sigma$, it follows that $mnG$ is open in $G$. This shows that $nH$ is open in $H$ and so is every subgroup of $H$ with index $n$.
\end{proof}


\begin{lemma}\label{erasebar}
Suppose that every finitely generated subgroup of $\A(K)$ has the congruence subgroup property. Then for any finitely generated subgroup $H$ of $\A(K)$ and any subset $J$ of $\A(K)$, we have $J\cap\overline{H}=J\cap H$.
\end{lemma}
\begin{proof}
Choose a finite subset $S$ of $\Omega_K$ such that $H\leq \A(O_S)$ and $S\cup\Omega=\Omega_K$.
Since $\A(O_S)$ has the congruence subgroup
property, $H$ is closed in $\A(O_S)$. Consequently,
$\A(O_S)\cap\overline{H}=H$. Now, since $\A(O_S)=\A(K)\cap \bigcap_{v\in \Omega}\A(O_v)$, while $\A(O_v)$ is closed in $\A(K_v)$,
$A(O_S)$ is closed in $\A(K)$, and hence
$\A(K)\cap\overline{\A(O_S)}=\A(O_S)$. This implies
$J\cap\overline{H}=J\cap\A(K)\cap\overline{\A(O_S)}\cap\overline{H}
=J\cap\A(O_S)\cap\overline{H}=J\cap H$.
\end{proof}

For any subgroup $J$ of $G:=\prod_{v\in\Omega}\A(K_v)$, its topological closure $\overline{J}$ is also a subgroup. In fact, since the map $q : G \times G \rarr G$ defined by $(P, Q) \mapsto P-Q$ is continuous, the preimage $q^{-1}(\overline{J})$ is a closed subset of $G \times G$ containing $J \times J$, and thus contains $\overline{J\times J}=\overline{J}\times \overline{J}$; then $q(\overline{J},\overline{J})=\overline{J}$ as desired. The following result generalizes Lemma 3.6 in \cite{PoonenVoloch}.

\begin{lemma}\label{tor}
If a finitely generated subgroup $H$ of $\A(K)$ has the congruence subgroup property, then every torsion element of $\overline{H}$ lies in $H$.
\end{lemma}
\begin{proof}
Write $H=T+F$, where $T$ is finite subgroup and $F$ is torsion-free. Suppose $a\in\overline{H}\setminus T$ with $ma=0$ for some nonzero integer $m$. Since $T$ is finite, there exists an open subgroup $U$ of $\prod_{v\in\Omega}\A(K_v)$ such that $(T+U)\cap (a+U)=\emptyset$.
Since $F$ is of finite index in $H$, by the congruence subgroup property of $H$, we may assume that $U\cap H\subset F$.
Lemma \ref{CSPsub} says that $U\cap H$ also has the congruence subgroup property, thus there exists an open subgroup $V$ of $\prod_{v\in\Omega}\A(K_v)$ such that $V\cap U\cap H=m(U\cap H)$. Because $ma=0$, the continuity of the multiplication-by-$m$ map ensures the existence of an open subgroup $W$ of $U$ such that $m(a+W)\subset U\cap V$. Now since $a\in\overline{H}$, there exist $t\in T$ and $f\in F$ such that $t+f\in H\cap (a+W)$. Consequently, $m(t+f)\in m(a+W)\cap H\subset U\cap V\cap H=m(U\cap H)$, whence
$m(t+f)=mf'$ for some $f'\in U\cap H\subset F$. Then $mt\in T\cap F=\{0\}$ and $m(f-f')=0$, which implies $f-f'\in T\cap F=\{0\}$ and $f=f'\in U$. This says $t+f\in (T+U)\cap (a+W)\subset (T+U)\cap (a+U)$, which is impossible.
\end{proof}

Suppose $J\subset \A(K)$ is a subgroup containing $H$. Then the inclusion $J\rarr\overline{J}$
canonically induces a group homomorphism
\begin{equation}\label{isoeq}
J/H\rarr\overline{J}/\overline{H}.
\end{equation}

\begin{lemma}\label{iso}
Suppose every finitely generated subgroup of $\A(K)$ has the congruence subgroup property.
Let $H\leq J$ be subgroups of $\A(K)$. If $H$ is
finitely generated, then {\em{(\ref{isoeq})}} is injective. If furthermore
the index $[J:H]$ is finite, then {\em{(\ref{isoeq})}} is
actually an isomorphism.
\end{lemma}
\begin{proof}
The first assertion follows from Lemma \ref{erasebar}. The congruence subgroup property of $J$ implies that
if $[J:H]$ is finite, then $H$ is open in $J$. Thus, $H=U\cap J$, for some open subgroup $U$ of $\prod_{v\in\Omega}\A(K_v)$.
Let $y$ be an arbitrary point of $\overline{J}$. Then $y\in z+U$ for some $z\in J$, and for any open subgroup $V$ of $U$,
$ (y-z+V)\cap J\neq \emptyset$. On the other hand, we have
$$
(y-z+V)\cap J\subset (y-z+V)\cap\left(U\cap J\right)=(y-z+V)\cap H.
$$

As the topology of $\prod_{v\in\Omega}\A(K_v)$ is generated by subgroups, it follows that $y-z\in\overline{H}$.
This shows the surjectivity of (\ref{isoeq}).
\end{proof}

The proof of the following proposition is postponed to Subsection \ref{pfCSPmain}.

\begin{prop}\label{CSPmain}
If $\A$ is isogenous to $\Gm^N\times \B$ for some nonnegative integer $N$ and some abelian variety $\B$ defined over $K$, then every finitely generated
subgroup of $\A(K)$ has the congruence subgroup property. In particular, the same conclusion holds if $\A$ is isogenous to an isotrivial semi-abelian variety defined over $K$.
\end{prop}

\subsection{The proof}\label{pfbase}

{\it Proof of Theorem \ref{AdelicIsotrivial_main} in the case where $\dim\X=0$}:

We write $\X=\Z$ to reflect this zero-dimensional situation.
As $\Z$ is zero-dimensional, by Lemma \ref{1stred}, we may replace $K$ by a finite extension if necessary and assume that every point of $\Z$ actually belongs to $\Z(K)$. In particular, the restriction $i_v|_{Z(K)}$ of the natural map $\xymatrix{A(K) \ar[r]^{i_v}& A(K_v)}$ is a bijection.

In view of Lemma \ref{erasebar}, we only have to show
$\prod_{v\in\Omega}\Z(K_v)\cap\overline{H}\subset\Z(K)$.
Let $J$ be the subgroup of $\A(K)$ generated by $H$ and $\Z(K)$.
By Proposition \ref{frob}, there exists a collection $\{U_n:n\geq 1\}$ of subgroups of $J$, which are open in every $v$-adic subtopology, such that $\bigcap_{n\geq 1} U_n$ is contained in the torsion subgroup of $J$. Let
$Q=(Q_v)_{v\in\Omega}$, with each $Q_v\in\Z(K_v)$, denote an element of $\prod_{v\in\Omega}\Z(K_v)$.
Suppose $Q$ is also contained in $\overline{H}$. Then there is a sequence $\{P_n\}_{n\geq 1}$ in $H$ such that at each $v$,
the sequence $\{i_v(P_n)\}_{n\geq 1}$ has $Q_v$ as its limit point
in $\A(K_v)$.  Write $Q_{(v)}=i_v^{-1}(Q_v)\in Z(K)$. Since each $U_n$ is open in the $v$-adic subtopology, for each $r\geq 1$, there exists an $N$ such that $P_n-Q_{(v)}\in U_r$, for $n\geq N$. It follows that for every pair $v,w\in\Omega$, the difference
$Q_{(v)}-Q_{(w)}$ belongs to $\bigcap_{r\geq 1} U_r$, whence a torsion point. Since the set $\{Q_{(v)}\}_{v\in\Omega}$ is finite, there exists a non-zero integer $m$ such that
$m(Q_{(v)}-Q_{(w)})=0$, for each pair $v$ and $w$.  Fix a $w\in\Omega$.
Then the difference $Q-Q_{(w)}=(Q_v-Q_w)_{v\in\Omega}\in\overline{J}$ is torsion, and hence, by Proposition \ref{CSPmain} and Lemma \ref{tor}, it is actually contained in $J$. In particular, $Q\in J$ is a global point.\qed

\begin{exm}\label{exm0}
The conclusion in Theorem \ref{AdelicIsotrivial_main} would fail, even in the case where $\dim\X=0$, if the hypothesis that
$H$ is finitely generated were removed. To see this, let $K$ be the field $\mathbb{F}_p(t)$ of rational functions over  $\mathbb{F}_p$. Fix a place $v_0$ of $K$
such that $t\notin O_{v_0}$. Let $\alpha, \beta\in K^*$ with $\beta-\alpha=\frac{a}{b}$, where $a,b\in \mathbb{F}_p[t]$. Denote by $\Z$ the  $K$-subvariety $\{\alpha, \beta\}$ of $\Gm$, and by $\xymatrix{K^* \ar[r]^{i_v}& K_v^*}$ the natural inclusion. Consider the sequence
\begin{equation}\label{xn}
x_n=\frac{(\prod_{i=1}^n \pi_i)^n+a}
{(\prod_{i=1}^n\pi_i)^{2n}+b}+\alpha,
\end{equation}
where $\pi_1, \pi_2, \pi_3,\ldots$ are all
irreducibles in $\mathbb{F}_p[t]$. Then the sequence $(x_n)$ has a limit $Q=(Q_v)_{v\in\Omega_K}$ in $\prod_{v\in\Omega_K}K_v^*$, where $Q_{v_0}=i_{v_0}(\alpha)$ and   $Q_v=i_v(\beta)$ for every $v\in\Omega_K\setminus\{v_0\}$, because
$$
x_n-\beta=\frac{\left(\prod_{i=1}^n \pi_i^n\right)\left(b-a\prod_{i=1}^n \pi_i^n\right)}
{b\left(\prod_{i=1}^n \pi_i^{2n}+b\right)}.
$$
\begin{enumerate}
\renewcommand{\labelenumi}{(\roman{enumi})}
\item If $\alpha\neq\beta$, then $Q\in\prod_{v\in\Omega_K}\Z(K_v)\cap\overline{\Gm(K)}\setminus\Z(K)$.
\item Suppose $\alpha=\beta\neq 1$ and set $b=1$ in (\ref{xn}). Then $\{x_n: n\geq 1\}\not\subset O_S^*$ for any finite $S\subset\Omega_K$; hence, by taking a subsequence, we may assume that every nonzero power of $x_n$ does not belong to the subgroup of $\Gm(K)$ generated by $\{\alpha,x_1,\ldots,x_{n-1}\}$. Letting $H$ be the subgroup of $\Gm(K)$ generated by $\{x_n: n\geq 1\}$,  we have $Q\in\prod_{v\in\Omega_K}\Z(K_v)\cap\overline{H}=\Z(K)\cap\overline{H}\setminus H$.
\end{enumerate}
\end{exm}

\section{The Inductive Step}\label{ind}
In this section, we complete the proof of Theorem \ref{AdelicIsotrivial_main} by reducing the general case to the zero-dimensional case which is established in Subsection \ref{pfbase}. Focusing on the case where the isogeny $f$ is the identity map until the very end of the reduction, the following Lemma \ref{int} and Lemma \ref{gen} are crucial to our inductive procedure. Their proofs, being long and independent of the rest materials in this section, will be postponed until
Subsection \ref{pfint}.

\begin{lemma}\label{int}
Let $N$ be a non-negative integer and $m$ a natural number. For each $v\in\Omega$, let $I_v$ be an ideal of
$K_v[X_0,\ldots,X_N]$, generated by elements of
$K_v^{p^m}[X_0,\ldots,X_N]$. Then $\bigcap_{v\in\Omega}
\left(I_v\cap K[X_0,\ldots,X_N]\right)$ is generated by elements of
$K^{p^m}[X_0,\ldots,X_N]$.
\end{lemma}

\begin{lemma}\label{gen} Let $N,m$ be non-negative integers. For each $v\in\Omega_K$, the ideal generated by those homogeneous polynomials in $K_v[X_0,\ldots,X_N]$
vanishing on a subset of $\Proj^N(K_v^{p^m})$ is actually generated by elements in $K_v^{p^m}[X_0,\ldots,X_N]$.
\end{lemma}

Then applications of the above are in order.

\begin{prop}\label{key} Let $m$ be a natural number. Let $\A_0$ be a semi-abelian variety defined over the largest finite subfield $\mathbb{F}_q\subset K$,
and $\X$  a closed $K$-subvariety which is not $K^{p^m}$-rational in $(A_0)_{\mathbb{F}_q}$.
Then there is a proper closed $K$-subvariety $\Y$ of $\X$
such that $\X(K_v)\cap \A_0(K_v^{p^m})\subset\Y(K_v)$ for all $v\in\Omega$.
\end{prop}
\begin{proof}
Since $\X$  is not $K^{p^m}$-rational in $(A_0)_{\mathbb{F}_q}$, we have an embedding of $\A_0$ into some $\Proj^N$ so that its underlying variety is
\begin{equation}\label{qvar1}
\left\{P\in\Proj^N(\overline{K}):
\begin{array}{l}
f(P)=0\text{ for all }f\in I\\ g(P)\neq 0\text{ for some }g\in J
\end{array}
\right\}
\end{equation}
for some homogeneous ideals $I$ and $J$ in $\mathbb{F}_q[X_0,\ldots,X_N]$, and that $\X$ is defined by (\ref{qvar1}) with $I$ replaced by a homogeneous radical ideal $I_X$ in $\overline{K}[X_0,\ldots,X_N]$ generated by elements of $K[X_0,\ldots,X_N]$, but not by those of $K^{p^m}[X_0,\ldots,X_N]$. For each $v\in\Omega$, consider the ideal $\tilde{I}_v$ in $K_v[X_0,\ldots,X_N]$ generated by homogeneous polynomials vanishing on the subset $\X(K_v)\cap\A_0(K_v^{p^m})$ of $\Proj^N(K_v)$. Let $\Y$ be the closed subvariety of $\A_0$ given by (\ref{qvar1}) except $I$ is replaced by
the homogeneous ideal
$$I_Y:=\left(\bigcap_{v\in\Omega} (\tilde{I}_v\cap K[X_0,\ldots,X_N])\right)\cap K^{p^m}[X_0,\ldots,X_N].$$
Thus  $\X(K_v)\cap \A_0(K_v^{p^n})\subset\Y(K_v)$ for all $v\in\Omega$.
We shall show $\Y\subsetneq \X$ by showing
\begin{equation}\label{temp}
I_X\subsetneq \overline{K}[X_0,\ldots,X_N]\cdot I_Y.
\end{equation}
To do so, we first apply Lemma \ref{gen} to deduce that $\tilde{I}_v$ is generated by elements in $K_v^{p^m}[X_0,\ldots,X_N]$.
Then by Lemma \ref{int}, we conclude that
$\bigcap_{v\in\Omega} (\tilde{I}_v\cap K[X_0,\ldots,X_N])$ is generated by elements in $K^{p^m}[X_0,\ldots,X_N]$, and that the right side of (\ref{temp}) equals to $\overline{K}[X_0,\ldots,X_N]\cdot\bigcap_{v\in\Omega} (\tilde{I}_v\cap K[X_0,\ldots,X_N])$. Since $I_X$ is not generated by elements of $K^{p^m}[X_0,\ldots,X_N]$, it proves (\ref{temp}).
\end{proof}

\begin{prop}\label{induct}
Let $\A_0$ be a semi-abelian variety defined over the largest finite subfield $\mathbb{F}_q\subset K$,
and $\X$ be a positive-dimensional closed $K$-subvariety of $\A_0$ such that all the largest dimensional
irreducible components of the translates $\X+P$, $P\in\A_0(\overline{K})$, are not $\overline{\mathbb{F}_p}$-rational in $(A_0)_{\mathbb{F}_q}$.
Let $H$ be a finitely generated subgroup of $\A_0(K)$,  Then there exists a closed $K$-subvariety $\Y$ of $\X$ with a
smaller dimension, satisfying $\prod_{v\in\Omega}\X(K_v)\cap
\overline{H}\subset\prod_{v\in\Omega}\Y(K_v)$.
\end{prop}
\begin{proof}
Let $\Frob:\A_0\rarr\A_0$ be the Frobenius endomorphism.
By taking  $H_0=\A_0(O_S)$ for some large enough finite $S\subset \Omega_K$, we assert that
that there is a finitely generated subgroup $H_0$ of $\A_0(K)$
such that $H\leq H_0$ and $\Frob(H_0)\leq H_0$. Since $\prod_{v\in\Omega}\X(K_v)\cap \overline{H}\subset\prod_{v\in\Omega}\X(K_v)\cap
\overline{H_0}$, it is enough to prove the desired result under the
additional hypothesis that $\Frob(H)\leq H$.

Assume that $\X$ is irreducible. First we apply the argument in the
proof of Theorem A. Part 1 in \cite{AV92} using the Hilbert scheme associated to equivalent compactification of $\A_0$, and conclude that there is a
positive integer $N$ such that for every $\gamma\in H$ the translate
$\X_{\gamma}=\X-\gamma$ is not $K^{p^N}$-rational in $(A_0)_{\mathbb{F}_q}$.  Therefore,
Proposition \ref{key} implies that there is a proper closed $K$-subvariety
$\Y_{\gamma}$ of $\X_{\gamma}$ such that $\X_{\gamma}(K_v)\cap
\A(K_v^{p^N})\subset\Y_{\gamma}(K_v)$ for all $v\in\Omega$.

Since the Frobenius endomorphism gives rise to an injection $\xymatrix{H\ar[r]^{\Frob^N} & H}$,
the index $[H:\Frob^N(H)]$ is finite, and hence Lemma \ref{iso} implies that there
are finitely many $\alpha_i$'s in $H$ such that
$\overline{H}=\bigcup_i\left(\alpha_i+\overline{\Frob^N(H)}\right)$.
Now, $\Y_{\alpha_i}+\alpha_i$ is a proper closed
$K$-subvariety of $\X$ such that
$\X(K_v)\cap\left(\alpha_i+\A_0(K_v^{p^N})\right)\subset\left(\Y_{\alpha_i}+\alpha_i\right)(K_v)$,
for all $v\in\Omega$. Then we prove the proposition by taking $Y=\bigcup_i\left(\Y_{\alpha_i}+\alpha_i\right)$.

In general, write $\X=\X_1\cup\dots\cup\X_m\cup \dots\cup\X_{m+n}$ where each $\X_i$ is irreducible and
$\mathrm{dim} \X_j=\mathrm{dim} \X$, for $j=1,...,m$; $\mathrm{dim} \X_i<\mathrm{dim} \X$, for $i=m+1,...,m+n$.
Then, for $j=1,...,m$, choose a closed proper $K$-subvariety $\Y_j$ of $\X_j$ satisfying
$\X_j(K_v)\cap \overline{H} \subset\Y_j(K_v)$ for all $v\in\Omega$. For $i=m+1,...,m+n$,
simply put $\Y_i=\X_i$. Then we complete the proof by taking $Y=\bigcup_i^{m+n}Y_i$.
\end{proof}

{\it Proof of Theorem \ref{AdelicIsotrivial_main}.} In view of Lemma \ref{1stred}, we may assume that the  isogeny $f:\A\rarr\A_0$ is defined over $K$, that $\A_0$ is defined over some finite subfield of $K$, and that every point in the kernel of $f$ lies in $\A(K)$.
If $\mathrm{dim} \X=0$, then the theorem is proved in Subsection \ref{pfbase}. In general, we prove by the induction on $\mathrm{dim} \X$.

Write $\X_0=f(\X)$. Proposition \ref{induct} applied to $\A_0$
ensures the existence of a closed $K$-subvariety $\Y_0$ of $\X_0$ of
smaller dimension such that
$$
\prod_{v\in\Omega}\X_0(K_v)\cap
\overline{f(H)}\subset\prod_{v\in\Omega}\Y_0(K_v).
$$
Write $\Y=f^{-1}(\Y_0)\cap \X$. Then the above implies
$$\prod_{v\in\Omega}\X(K_v)\cap
\overline{H}\subset \prod_{v\in\Omega}\X(K_v)\cap \prod_{v\in\Omega}f^{-1}(\Y_0)(K_v)=\prod_{v\in\Omega}\Y(K_v).$$
The assumption in Theorem \ref{AdelicIsotrivial_main} is preserved when $\X$ is replaced by $\Y$, hence
the induction hypothesis implies that $Y(K)\cap H$ is finite and
$$\prod_{v\in\Omega}\Y(K_v)\cap
\overline{H}=Y(K)\cap H.$$
Therefore,
$$\prod_{v\in\Omega}\X(K_v)\cap
\overline{H}\subset \Y(K)\cap H\subset \X(K)\cap H\subset \prod_{v\in\Omega}\X(K_v)\cap
\overline{H} .$$
This completes the proof.
\qed

In order to deduce Corollary \ref{ss2curve} from Theorem \ref{AdelicIsotrivial_main}, we need the following result.
\begin{lemma}\label{sandwich}
Let $C_1\rightarrow C\rightarrow C_0$ be a chain of nonconstant maps between projective curves defined over $\overline{K}$ with $C$ smooth. Suppose that both $C_0$ and $C_1$ as well as the composition $C_1\rightarrow C_0$ are defined over $\overline{\mathbb{F}_p}$.
Then $C$ is also defined over $\overline{\mathbb{F}_p}$.
\end{lemma}
\begin{proof}
The given chain of maps induces the following diagram of their function fields:
\begin{center}
\begin{picture}(100,140)(20,-10)
\put(0,120){$\overline{\mathbb{F}_p}(C_1)$} \put(30,123){\vector(1,0){65}} \put(100,120){$\overline{K}(C_1)$}
                            \put(100,60){$\overline{K}(C)$}
                            \put(105,75){\line(0,1){40}}
\put(5,15){\line(0,1){100}} \put(105,15){\line(0,1){40}}
\put(0,0){$\overline{\mathbb{F}_p}(C_0)$} \put(30,3){\vector(1,0){65}} \put(100,0){$\overline{K}(C_0)$,}
\end{picture}
\end{center}
where both columns are finite extensions, and the maps in both rows are $\otimes_{\overline{\mathbb{F}_p}} \overline{K}$.
To prove this lemma, since $C$ is smooth, it suffices to find a field $F$ with transcendence degree $1$ over $\overline{\mathbb{F}_p}$ such that $F\otimes_{\overline{\mathbb{F}_p}} \overline{K}$ is $\overline{K}$-isomorphic to $\overline{K}(C)$. First we assume that $\overline{\mathbb{F}_p}(C_1)$ is separable over $\overline{\mathbb{F}_p}(C_0)$. Let $N$ be the normal closure of $\overline{\mathbb{F}_p}(C_1)$ over $\overline{\mathbb{F}_p}(C_0)$. Identifying all fields involved as subfields of $N\otimes_{\overline{\mathbb{F}_p}} \overline{K}$, we take $F=N\cap \overline{K}(C)$. Galois theory shows that $[N\otimes_{\overline{\mathbb{F}_p}} \overline{K}: \overline{K}(C)]=[N:F]$. Since $F\otimes_{\overline{\mathbb{F}_p}} \overline{K}\subset \overline{K}(C)$ and $[N\otimes_{\overline{\mathbb{F}_p}} \overline{K}: F\otimes_{\overline{\mathbb{F}_p}} \overline{K}]\leq [N: F]$, we conclude that $F\otimes_{\overline{\mathbb{F}_p}} \overline{K}=\overline{K}(C)$ as desired.

In the general case, let $L$ be the separable closure of $\overline{\mathbb{F}_p}(C_0)$ in $\overline{\mathbb{F}_p}(C_1)$.
The preceding argument yields a field $F'$ with transcendence degree $1$ over $\overline{\mathbb{F}_p}$ such that $F'\otimes_{\overline{\mathbb{F}_p}} \overline{K}$ is $\overline{K}$-isomorphic to $\overline{K}(C)\cap\left(L\otimes_{\overline{\mathbb{F}_p}} \overline{K}\right)$. Since $\overline{K}(C)$ is purely inseparable over $\overline{K}(C)\cap\left(L\otimes_{\overline{\mathbb{F}_p}} \overline{K}\right)$, the property of the Frobenius map shows that $\overline{K}(C)\cap\left(L\otimes_{\overline{\mathbb{F}_p}} \overline{K}\right)=\overline{K}(C)^q$ for some power $q$ of $p$, and since $\overline{\mathbb{F}_p}$ is perfect, the field $F=F'^{\frac{1}{q}}$ is the one we look for.
\end{proof}

{\it Proof of Corollary \ref{ss2curve}.}
Let $X$ be a non-isotrivial smooth projective $K$-curve with its Jacobian $J$ isogenous to an isotrivial abelian variety $\A_0$.
Without loss of generality, we can assume that $A_0$ is defined over $\overline{\mathbb{F}_p}$.
Denote by $f:J\rarr\A_0$ the isogeny and by $\breve{f}:\A_0\rarr J$ its dual. Let $m$ be the positive integer such that $f\circ\breve{f}$ is the multiplication-by-$m$ map on $\A_0$, and $\breve{f}\circ f$ is the multiplication-by-$m$ map on $J$.
In view of the discussion given in Section \ref{AdelicIsotrivial_intro} , we need to show that $\Sha(J)$ is finite and each translate
$f(X)+P$, $P\in A_0(\overline{K})$ is not $\overline{\mathbb{F}_p}$-rational.

Now, since $\A_0$ is isotrivial, $\Sha(A_0)$ is finite, by  Tate \cite{Sha}.
Choose a prime $l$ not dividing $pm$. The isogenies $f$ and $\breve{f}$ induce a chain  $$\Sha(J)[l^{\infty}]\rarr\Sha(\A_0)[l^{\infty}]\rarr\Sha(J)[l^{\infty}]$$
of maps between the $l$-primary part of $\Sha(J)$ and of $\Sha(\A_0)$ such that the composition is an isomorphism. In particular, $\Sha(\A_0)[l^{\infty}]\rarr\Sha(J)[l^{\infty}]$ is surjective, hence $\Sha(J)[l^{\infty}]$ is finite as $\Sha(\A_0)[l^{\infty}]$ is. By another result of Tate \cite{Sha}, it follows that $\Sha(J)$ is finite as desired.

Suppose $C_0:=f(X)+P$ is $\overline{\mathbb{F}_p}$-rational in $\A_0$. Write $C=X+Q$, for some $Q\in f^{-1}(P)$ and
let $C_1$ be an irreducible component of the pre-image of $C_0$ under $\xymatrix{\A_0\ar[r]^m & \A_0}$. Then Lemma \ref{sandwich}
is applicable to the chain
$\xymatrix{C_1\ar[r]^{\breve{f}} & C \ar[r]^f & C_0}$. Consequently, $C$ is defined over $\overline{\mathbb{F}_p}$,
and hence $X$, being isomorphic to $C$, is isotrivial. This is a contradiction.\qed

\begin{exm}\label{exm1}
The conclusion in Theorem
\ref{AdelicIsotrivial_main} would fail if no hypothesis were put on $\X$. To see this, let $K$ be the field $\mathbb{F}_p(t)$ of rational functions over  $\mathbb{F}_p$, and $H=\la t\ra$ be the cyclic subgroup of $\Gm(K)$
generated by $t$. Take a cofinite subset $\Omega$ of $\Omega_K$ such
that $t\in O_v^*$ for every $v\in\Omega$. For any
$m\geq n$, we have
$$
t^{p^{m!}}-t^{p^{n!}}=\left(t^{p^{n!(\frac{m!}{n!}-1)}}-t\right)^{p^{n!}}.
$$
Thus, the sequence
$(t^{p^{n!}})_{n\geq 1}$ in $H$ is Cauchy, and admits a limit $Q=(Q_v)_{v\in\Omega}$ in $\overline{H}$ by compactness. Note that $Q_{v_{t-1}}=1$, where $v_{t-1}\in\Omega$ is the unique one satisfying $t-1\in m_{v_{t-1}}$; while $Q_v\neq 1$ for each $v\in\Omega\setminus\{v_{t-1}\}$. Hence $Q\in\prod_{v\in\Omega}\Gm(K_v)\cap
\overline{H}\setminus\Gm(K)$.
\end{exm}

\section{The proofs of key intermediate results}\label{pfkey}
\subsection{The proof of Proposition \ref{CSPmain}}\label{pfCSPmain}
In this subsection, we fix a finitely generated subgroup $H\subset A(K)$.
The number field counter part of the following lemma (for $\Omega$ consisting of only non-Archimedean places) is just a reinterpretation of Theorem 1, \cite{Chev}, and it can actually be carried over to the function field case. I thank the referee for pointing out the present much shorter proof using Galois cohomology.

\begin{lemma}\label{chev}
If $\A=\Gm$, then every subgroup of $H$ of index prime to $p$ is open.
\end{lemma}
\begin{proof}
In view of Lemma \ref{CSPsub}, we only need to consider the case where $H=\Gm(O_T)=O_T^*$ for a finite $T\subset \Omega_{K}$. For any finite subset $S\subset\Omega$, consider the open subgroup
$U_S=\prod_{v\in S}1+m_v$ of $\prod_{v\in S}K_v^*$. We shall
proves the lemma by showing that for any natural number $m$ prime to $p$, there is some $U_S$ such that $O_T^*\cap U_S\subset (O_T^*)^m$.

Now, Kummer theory gives rise to the following commutative diagram
$$
\xymatrix{O_T^*/(O_T^*)^m \ar@{^{(}->}[r] \ar[d] & K^*/(K^*)^m \ar[r]^{\sim} \ar[d] & H^1(K,\mu_m)\ar[d]\\
\displaystyle{\prod_{v\in\Omega}} O_v^*/(O_v^*)^m \ar@{^{(}->}[r] & \displaystyle{\prod_{v\in\Omega}} K_v^*/(K_v^*)^m \ar[r]^{\sim} & \displaystyle{\prod_{v\in\Omega}} H^1(K_v,\mu_m),}\\
$$
where the two injections are clear. As the Galois group of $K(\mu_m)/K$ is cyclic, by Lemma I.9.3 of \cite{ADT}, the right vertical arrow is an injection, and hence so is the left one. Since
$O_T^*/(O_T^*)^m$ is finite, there exists a finite subset $S\subset\Omega$  such that the left vertical arrow induces an injection $O_T^*/(O_T^*)^m \hookrightarrow\prod_{v\in S} O_v^*/(O_v^*)^m$. Hensel's lemma shows
$U_S
\subset \prod_{v\in S}(O_v^*)^m$, whence $O_T^*\cap U_S\subset (O_T^*)^m$ as desired.
\end{proof}

\begin{lemma}\label{ppower}
If $\A=\Gm$, then very subgroup of $H$ of  $p$-power index is open in the $v$-adic subtopology, for every $v\in\Omega_K$.
\end{lemma}
\begin{proof}
Again, by Lemma \ref{CSPsub}, we only need to consider the case where $H=O_T^*$.
Now, we have $(O_T^*)^{p^e}=O_T^* \cap (K^*)^{p^e}$, which is shown in the proof of Proposition \ref{frob}
to be an open subgroup in $O_T^*$ in the $v$-adic subtopology, for every $v\in\Omega_K$.
\end{proof}

\begin{cor}\label{GmCSP}
If $\A=\Gm^N$, then every finitely generated subgroup of $\A(K)$ has the congruence subgroup property.\qed
\end{cor}

In the case where $\A$ is an abelian variety, the first assertion of Proposition \ref{CSPmain} is essentially proved by Milne, who generalizes a result of Serre \cite{SerreCong} in the case where $K$ was a number field.

\begin{prop}\label{ACSP}
Suppose that $\A$ is an abelian variety defined over $K$. Then every  subgroup $H$ of $\A(K)$ has the congruence subgroup property.
\end{prop}
\begin{proof}
The case where $H=\A(K)$ is exactly Corollary 1, \cite{Mil}. Then other cases follow from Lemma \ref{CSPsub}.
\end{proof}

{\it Proof of Proposition \ref{CSPmain}}: Let $H$ be a finitely generated subgroup of
$\A(K)$. We need to show that $nH\subset H$ is open for every $n$.
Suppose we are given the isogeny $\phi:\A\rarr \Gm^N\times \B$. It follows from Corollary \ref{GmCSP} and Proposition \ref{ACSP} that $\phi(H)$ has the congruence subgroup property. In particular $mn\phi(H)\subset \phi(H)$ is an open subgroup, for every $m$.
Denote
$$T=\text{ker}\left[\xymatrix{\A(K)\ar[r]^-{\phi} & \Gm^N(K)\times \B(K) }  \right].$$
Since $T$ is finite, Lemma \ref{adehaus} implies the existence of an open subgroup $U$ of $H$ of a
finite index $m$ such that $U\cap T$ is trivial.

Now, since $nmH+T=\phi^{-1}(nm\phi(H))$ is open in $H+T=\phi^{-1}(\phi(H))$, the subgroup
$H\cap\left(nmH+T\right)$ is open in $H$. Also, since $nmH\subset mH\subset U$ and $U\cap T$ is trivial, we see that $U\cap\left(nmH+T\right)=nmH$. Therefore $nmH$ is open in $H$, and thus
so is $nH$.
\qed

\subsection{The proofs of Lemma \ref{int} and Lemma \ref{gen}}\label{pfint}
Our tool for proving Lemma \ref{int} and Lemma \ref{gen} is the iterative derivation. An iterative derivation on a field $L$ is a sequence $\{D^{(i)}\}_{i\geq 0}$ of elements in the
$L$-algebra of additive endomorphisms on $L$
such that
\begin{enumerate}
\renewcommand{\labelenumi}{(\theenumi)}
\renewcommand{\theenumi}{\roman{enumi}}
\item $D^{(0)}$ is the identity operator.
\item $D^{(i)}(xy)=\sum_{j=0}^i D^{(j)}(x)D^{(i-j)}(y)$, for $i\geq 0$ and $x,y\in L$.
\item\label{iter} $D^{(i)}D^{(j)}= {i+j\choose i}D^{(i+j)}$ for $i,j\geq 0$, where $D^{(i)}D^{(j)}$ denotes the composition of $D^{(i)}$ and $D^{(j)}$, and the rational integer ${i+j \choose i}$ is the binomial coefficient.
\end{enumerate}

Assume that $L$ is of characteristic $p$. Then the following Lucas's lemma (see, for example, \cite{SFK}),
is useful for telling if ${i+j \choose i}\not=0$ in $L$. For each nonnegative integer $i$, let $i=\sum_{n=0}^d i_n p^n$, $0\leq i_n< p$, denote
its base $p$ expansion.
\begin{lemma}\label{bino}
The binomial coefficient $i\choose j$ is not divisible by $p$ if and
only if $i_n\geq j_n$ for all $n$.\qed
\end{lemma}

The defining property (\ref{iter}) implies $D^{(i)}\circ D^{(j)}=D^{(j)}\circ D^{(i)}$.
Also, repeated applications of the property (\ref{iter})
gives
\begin{equation*}
\prod_{n=0}^d \left(D^{(p^n)}\right)^{i_n}
=c_i D^{(i)},
\end{equation*}
where
\begin{equation*}
c_i=\prod_{n=0}^d\left[{\sum_{s=0}^n i_s p^s\choose i_n p^n}\prod_{a=1}^{i_n} {ap^n\choose p^n}\right].
\end{equation*}
Now, Lemma \ref{bino} implies $c_i\in L^*$, and hence
\begin{equation}\label{ci*}
D^{(i)}=c_i^{-1}\cdot \prod_{n=0}^d \left(D^{(p^n)}\right)^{i_n}.
\end{equation}
Inspired by the proof of Claim 2.2.3,
\cite{Ogus}, we consider the operator
$$\Delta_m:=\sum_{i=0}^{p^m-1} (-t)^i D^{(i)}$$
on $L$ for some $t\in L$ satisfying
\begin{equation}\label{t}
D^{(i)}\left((-t)^j\right)=(-1)^i {j\choose i}
t^{j-i}\qquad\qquad\text{for each }i,j\geq 0.
\end{equation}

For each $m\geq 0$, let $L_m=\{x\in L:D^{(l)}(x)=0, \;\text{ if }\;
1\leq l< p^m\}$,
which is a subfield of $L$.
\begin{lemma}\label{magicsum}
For every $c\in L$ and every $m\geq 0$, the element
$\Delta_m(c)\in L_m$.
\end{lemma}
\begin{proof}
In view of \eqref{ci*}, we only need to show that for every natural number $s<m$,
$$
D^{(p^s)}(\Delta_m(c))=0.
$$
For simplicity, set
$j=p^s$. It follows from the property (\ref{iter}) and the assumption \eqref{t}
that
$$
D^{(j)}(\Delta_m(c))=\sum_{i=0}^{p^m-1}\sum_{l=0}^j
(-1)^l {i\choose l}{i+j-l\choose i}(-t)^{i-l}
D^{(i+j-l)}(c).
$$
Lemma \ref{bino} implies that ${i\choose
l}{i+j-l\choose i}$ is a multiple of $p$ unless
both $i_n\geq l_n$ and $j_n\geq l_n$ hold for all
$n$, which occurs only when $l\in \{0,j\}$, since
$j_s=1$ and $j_n=0$ for all $n\neq s$. We also
note that in case
where $l=j$, those terms with $i<j$ vanish as
${i\choose j}=0$. Putting these together, we
obtain
\begin{equation*}
\begin{array}{rl}
  & D^{(j)}(\Delta_m(c)) \\
= & \sum_{i=0}^{p^m-1}
{i+j\choose i}(-t)^{i} D^{(i+j)}(c)+
\sum_{i=0}^{p^m-1} (-1)^j {i\choose j}(-t)^{i-j}
D^{(i)}(c) \\
= & \sum_{i=j}^{j+p^m-1} {i\choose
j}(-t)^{i-j} D^{(i)}(c)+ \sum_{i=j}^{p^m-1}
(-1)^j {i\choose j}(-t)^{i-j} D^{(i)}(c) \\
= & \sum_{i=p^m}^{p^m+p^s-1} {i\choose
p^s}(-t)^{i-p^s}D^{(i)}(c),
\end{array}
\end{equation*}
where the last equality holds because
$1+(-1)^j=1+(-1)^p=0$ in $K$. Finally,
since $s<m$, for every $i$ satisfying $p^m\leq i\leq p^m+p^s-1$, we have $i_s=0$,
and hence ${i\choose
p^s}$ is a multiple of $p$, by Lemma \ref{bino}. Thus, each term in
the last sum vanishes. This finishes the proof.
\end{proof}

For each $i\geq 0$, we extend $D^{(i)}$ to an additive endomorphism on the polynomial ring $L[X_0,\ldots,X_N]$ by sending $X_i$ to $0$
for every $i\in\{0,1,\ldots,N\}$. It is easy to
verify that for all $i\geq 0$, $f,g\in L[X_0,\ldots,X_N]$,
$$D^{(i)}(fg)=\sum_{j=0}^i D^{(j)}(f)D^{(i-j)}(g),,$$
and, for all $m\geq 0$,
$$L_m[X_0,\ldots,X_N]=\{g\in
L[X_0,\ldots,X_N]:D^{(l)}(g)=0,\;\text{ if } 1\leq l< p^m\}.$$

\begin{lemma}\label{der}
For any  positive integer $m$,  an
ideal $I$ of $L[X_0,\ldots,X_N]$ is
generated by elements of $L_m[X_0,\ldots,X_N]$
if and only if the condition $D^{(i)}(I) \subset
I$ holds for all $1\leq i< p^m$.
\end{lemma}
\begin{proof}

Suppose
$D^{(i)}(I) \subset I$ for all $1\leq i< p^m$. Let
$J$ be the ideal of $L[X_0,\ldots,X_N]$ generated
by $I \cap L_m[X_0,\ldots,X_N]$. To complete the proof, we only need to show $I=J$, as the implication in the opposite direction is clear.
Choose a lexicographic order on the set of
monomials in $X_0,\ldots,X_N$. With respect to this
order, for each non-zero polynomial $f\in K[X_0,\ldots,X_N]$, the \em degree \em of $f$ is defined to be the largest monomial appearing in the
expression of $f$ with a non-zero coefficient, and $f$ is  \em monic \em if this coefficient is $1$.

Suppose that
$I\setminus J$ is a non-empty set and let $f\in I\setminus J$ be an element of the smallest degree.
We also choose $f$ to be monic.
Since $D^{(i)}(1)=0$, for all positive integer $i$, the degree of $D^{(i)}(f)$ is smaller than that of $f$.
Consequently, $D^{(i)}(f) \notin I
\setminus J$, by the choice of $f$. On the other hand, since for every
$1\leq i< p^m$, $D^{(i)}(f)\in D^{(i)}(I)
\subset I$, we must have $D^{(i)} (f)\in J$.
Now, consider the element
$$g=f+ \sum_{i=1}^{p^m-1}
(-t)^i D^{(i)}(f)\in L[X_0,\ldots,X_N].$$
By the above argument, $g\in I$, and by Lemma
\ref{magicsum}, $g\in L_m[X_0,\ldots,X_N]$. Hence, $g\in J$ and
$f=g-\sum_{i=1}^{p^m-1}(-t)^i D^{(i)}(f)\in J$, a contradiction.
\end{proof}

Now we construct a desired iterative derivation on $K$.
Choose an element $t\in K$ such
that $K$ is a finite separable extension of the
function field $\mathbb{F}_p(t)$ of one variable over
$\mathbb{F}_p$. Choose a place $v_0\in\Omega_K$
which restricts to a place $w\in\Omega_{\mathbb{F}_p(t)}$ corresponding to a separable irreducible
polynomial in $\mathbb{F}_p[t]$ such that $\mathbb{F}_p(t)_w=K_{v_0}$. Let $\alpha$ be a
root of this  polynomial. Then $\mathbb{F}_p(t)_w$ is a natural subfield of $\mathbb{F}_p(\alpha)\left((t-\alpha)\right)$ and we have
a tower $\mathbb{F}_p(\alpha)(t)\subset
K(\alpha)\subset
\mathbb{F}_p(\alpha)\left((t-\alpha)\right)$ of
fields.
By Remark 1 in \cite{VolochWronskians}, there exists an iterative derivation $\{D_{K(\alpha)}^{(i)}\}_{i\geq 0}$ on
$K(\alpha)$ such that
$D_{K(\alpha)}^{(j)}\left((t-\alpha)^i\right)={i\choose
j}(t-\alpha)^{i-j}$ and
$\left(K(\alpha)\right)^{p^m}=\{x\in
K(\alpha):D_{K(\alpha)}^{(l)}(x)=0,\;\text{ if }\; 1\leq l< p^m\}$,
for $i,j,m\geq 0$. Denoting by $D_K^{(i)}$ the restriction of $D_{K(\alpha)}$ on $K$, we get an iterative derivation $\{D_K^{(i)}\}_{i\geq 0}$ on
$K$.
It is not hard to check
that $D_K^{(j)}(t^i)={i\choose j}(t)^{i-j}$, whence \eqref{t} holds for $D=D_K$. Also, from the
separability assumption, we have
$K^{p^m}=\{x\in K:D_K^{(l)}(x)=0\text{ if } 1\leq l<
p^m\}$, for $m\geq 0$. Moreover, using the fact $[K:K^p]=p$, one can show that for each $i$,  the endomorphism $D_K^{(i)}$ is continuous with respect to every
place of $K$. Therefore, for each place
$v\in\Omega$, we extend $\{D_K^{(i)}\}_{i\geq 0}$ and obtain an
iterative derivation $\{D_{K_v}^{(i)}\}_{i\geq 0}$
on $K_v$.

{\it Proof of Lemma \ref{int}}:
 Since $I_v$ is generated by
elements of $K_v^{p^m}[X_0,\ldots,X_N]$, which lie in the kernel of  those $D_{K_v}^{(i)}$ with $1\leq i
< p^m$, it follows that for these $i$ we have
$D_{K_v}^{(i)}(I_v)\subset I_v$ for each
$v\in\Omega$. But then
$$
\begin{array}{rl}
 & D_K^{(i)}\left(\bigcap_{v\in\Omega}
\left(I_v\cap
K[X_0,\ldots,X_N]\right)\right)\\
\subset &\bigcap_{v\in\Omega}
\left(D_{K_v}^{(i)}\left(I_v\cap
K[X_0,\ldots,X_N]\right)\right)\\
\subset &
\bigcap_{v\in\Omega}
\left(D_{K_v}^{(i)}(I_v)\cap
K[X_0,\ldots,X_N]\right)\\
\subset &\bigcap_{v\in\Omega}
\left(I_v\cap K[X_0,\ldots,X_N]\right)
\end{array}
$$
for all
$1\leq i< p^m$. Then we complete the proof by applying Lemma \ref{der}.\qed
\newline

{\it Proof of Lemma \ref{gen}}:
Fix a subset $\Sigma$ of $\Proj^N(K_v^{p^m})$ for some place $v\in\Omega_K$ and some positive integer $m$, and denote by $I_v$ the ideal in $K_v[X_0,\ldots,X_N]$ generated by homogeneous polynomials which vanish on $\Sigma$.
Let $f\in I_v$ be a homogeneous polynomial and $P\in \Sigma$. By the definition of $D_{K_v}^{(i)}(f)$ and the assumption $P\in\Proj^N(K_v^{p^m})$, we  have $0=D_{K_v}^{(i)}\left(f(P)\right)=D_{K_v}^{(i)}(f)\left(P\right)$ for each $1\leq i< p^m$. This shows $D_{K_v}^{(i)}(I_v)\subset I_v$ for all $1\leq i< p^m$. Again, we complete the proof by using Lemma \ref{der}.

\section*{acknowledgements}
I thank the referee for his/her careful comments on an earlier version of this paper, which generalizes a part of my dissertation completed under the advising of Jos\'e Felipe Voloch. Also, I appreciate that Andreas Schweizer, Ki-Seng Tan, and Julie Wang made many useful suggestions on my revision. Especially, I am grateful to Tan for showing me examples of succinct writing.

\bibliographystyle{alpha}
\bibliography{mybib}

\end{document}